\documentstyle[fleqn,12pt]{article}
\begin{document}\pagenumbering{arabic}\setcounter{page}{1}
\pagestyle{plain}
\baselineskip=16pt

\thispagestyle{empty}
\rightline{MSUMB 97-05, September 1997} 
\vspace{0.3cm}

\begin{center}
{\Large\bf 
Comment on the differential calculus on the quantum exterior plane }
\end{center}

\vspace{1cm}
\begin{center} Salih Celik$^1$, Sultan A. Celik$^{2,1}$ and Metin Arik$^3$ 
\end{center}

\noindent
$^1$ {\footnotesize Mimar Sinan University, Department of Mathematics, 
80690 Besiktas, Istanbul, TURKEY.}\\
$^2$ {\footnotesize Yildiz Technical University, Department of Mathematics, 
Sisli, Istanbul, TURKEY. }\\
$^3$ {\footnotesize Bogazici University, Department of Physics, Bebek, 
Istanbul, TURKEY. }

\vspace{1.5cm}
{\bf Abstract}

We give a two-parameter quantum deformation of the exterior plane and 
its differential calculus without the use of any R-matrix and relate 
it to the differential calculus with the R-matrix. We prove that there 
are two types of solutions of the Yang-Baxter equation whose symmetry 
group is $GL_{p,q}(2)$. We also give a two-parameter deformation of 
the fermionic oscillator algebra. 

\vfill\eject
Quantum groups are a generalization of the concept of groups. 
During the past few years, these new mathematical objects have 
found wide interest among theoretical physicsts and mathematicians. 
More precisely, the quantum group is an example of a Hopf algebra 
which is related to noncommutative geometry.$^1$ After Wess-Zumino 
introduced the differential calculus on the quantum (hyper)plane, 
the quantum plane was generalized to the supersymmetric quantum 
(super)plane by some authors.$^{2-5}$ 

In this letter we shall give a differential calculus on the quantum 
exterior plane whose symmetry group is $GL_{p,q}(2)$ and obtain 
two R-matrices. They are both solutions of the Yang-Baxter equation. 

Let us begin with the quantum exterior (dual) plane, defined as 
the polynomial ring generated by coordinates $\theta$, $\phi$ 
which satisfy$^6$ 
$$ \theta \phi + p^{-1} \phi \theta = 0 \eqno(1) $$
$$ \theta^2 = 0 = \phi^2 \eqno(2)$$
where $p$ is a complex deformation parameter. To develop the differential 
calculus on the $p$-exterior plane in terms of the coordinates satisfying 
(1), (2), we shall make the following ansatz for the commutation relations 
of the coordinates with their differentials. Let the differentials of 
coordinates be denoted by 
$$ \Theta = d\theta \qquad \Phi = d\phi. \eqno(3)$$
In Ref. 1, Wess and Zumino have interpreted $\theta$ and $\phi$ as the 
differentials of coordinates of the quantum plane, respectively. Recall 
that the quantum plane is defined as the polynomial ring generated by 
coordinates $x$, $y$ obeying the relation 
$$x y = q y x.$$
As an alternative to this interpration, $x$ and $y$ can be identified with 
the differentials of $\theta$ and $\phi$ as in Ref. 6. The definition (3) 
is closely related to this approach. We assume that 
$$\theta \Theta = A \Theta \theta \qquad 
  \theta \Phi = F_{11} \Phi \theta + F_{12} \Theta \phi $$
$$\phi \Phi = B \Phi \phi \qquad 
  \phi \Theta = F_{21} \Theta \phi + F_{22} \Phi \theta. \eqno(4)$$
Then we find, from (2), 
$$A = 1 \qquad B = 1. \eqno(5)$$
The consistency condition 
$$ d(\theta \phi + p^{-1} \phi \theta) = 0 $$
$$ (\theta \phi + p^{-1} \phi \theta)\Theta = 0 \qquad 
   (\theta \phi + p^{-1} \phi \theta)\Phi = 0 \eqno(6)$$
gives 
$$F_{22} + p F_{11} = 1 \qquad F_{21} + p F_{12} = p \qquad 
  F_{12} F_{22} = 0. \eqno(7)$$
We define the exterior derivative {\sf d} obeying the condition: 
$${\sf d}^2 = 0 \eqno(8)$$
and the graded Leibniz rule 
$${\sf d}(f g) = ({\sf d} f) g + (- 1)^{\hat{f}} f ({\sf d} g) \eqno(9)$$
where $\hat{f} = 0$ for even variables and $\hat{f} = 1$ for odd variables. 
Applying the exterior differential {\sf d} on second and fourth relations 
of eq.(4) and using (8), (9) one gets 
$$F_{11} F_{21} = 1 - F_{12} - F_{22}. \eqno(10)$$
We now assume that the commutation relation of differentials has the form 
$$\Theta \Phi = q \Phi \Theta \eqno(11)$$
where $q$ is another complex deformation parameter. Then we have
$$F_{11} = q (1 - F_{12}) \qquad F_{21} = q^{-1} (1 - F_{22}). \eqno(12)$$
The system (7), (10), (12) admits two solutions: 

{\it Type I} 
$$A = 1 \qquad F_{11} = q \qquad F_{21} = p $$
$$B = 1 \qquad F_{12} = 0 \qquad F_{22} = 1 - pq \eqno(13)$$
and the relations (4) take the form 
$$\theta \Theta = \Theta \theta \qquad 
  \theta \Phi = q \Phi \theta $$
$$\phi \Phi = \Phi \phi \qquad 
  \phi \Theta = p \Theta \phi + (1 - pq) \Phi \theta. \eqno(14)$$

{\it Type II} 
$$A = 1 \qquad F_{11} = p^{-1} \qquad F_{21} = q^{-1} $$
$$B = 1 \qquad F_{12} = 1 - p^{-1}q^{-1} \qquad F_{22} = 0 \eqno(15)$$
and in this case the relations (4) take the form 
$$\theta \Theta = \Theta \theta \qquad 
  \theta \Phi =  p^{-1} \Phi \theta - (1 - p^{-1} q^{-1} \Theta \phi$$
$$\phi \Phi = \Phi \phi \qquad 
  \phi \Theta = q^{-1} \Theta \phi. \eqno(16)$$

To complete the differential geometric scheme we introduce derivatives 
of the quantum exterior plane in the standard way 
$${\sf d} = \Theta \partial_\theta + \Phi \partial_\phi. \eqno(17)$$
Multiplying this expression from the right by $\theta f$ and $\phi f$, 
respectively, and using the graded Leibniz rule for partial derivatives 
$$\partial_i(f g) = (\partial_i f)g + (-1)^{\hat{f}} f (\partial_i g) 
  \eqno(18)$$
one finds 
$$ \partial_\theta \theta = 1 - \theta \partial_\theta - 
   F_{12} \phi \partial_\phi 
   \qquad \partial_\theta \phi = - F_{21} \phi \partial_\theta$$
$$ \partial_\phi \phi = 1 - \phi \partial_\phi - 
   F_{22} \theta \partial_\theta 
 \qquad \partial_\phi \theta = - F_{11} \theta \partial_\phi. \eqno(19)$$
These commutation relations, for type I, are 
$$ \partial_\theta \theta = 1 - \theta \partial_\theta 
   \qquad \partial_\theta \phi = - p \phi \partial_\theta$$
$$ \partial_\phi \phi = 1 - \phi \partial_\phi + 
   (pq - 1) \theta \partial_\theta 
   \qquad \partial_\phi \theta = - q \theta \partial_\phi \eqno(20\mbox{a})$$
and for type II, are 
$$ \partial_\theta \theta = 1 - \theta \partial_\theta + 
   (p^{-1} q^{-1} - 1) \phi \partial_\phi 
   \qquad \partial_\theta \phi = - q^{-1} \phi \partial_\theta$$
$$ \partial_\phi \phi = 1 - \phi \partial_\phi 
   \qquad \partial_\phi \theta = - p^{-1} \theta \partial_\phi. 
  \eqno(20\mbox{b})$$

The commutation relations between the derivatives can be easily obtained 
by using that ${\sf d}^2 = 0$. So it follows that 
$$0 = {\sf d}^2 = \Theta^2 \partial_\theta^2 + \Phi \Theta (\partial_\theta 
  \partial_\phi + q \partial_\phi \partial_\theta) + 
  \Phi^2 \partial_\phi^2 $$
which says that 
$$\partial_\theta \partial_\phi + q \partial_\phi \partial_\theta = 0 
  \qquad \partial_\theta^2 = 0 = \partial_\phi^2. \eqno(21)$$
Finally to find the commutation rules between the differentials 
and derivatives we shall assume that they have the following form 
$$ \partial_\theta \Theta = A_{11} \Theta \partial_\theta + 
   A_{12} \Phi \partial_\phi $$
$$ \partial_\theta \Phi = A_{21} \Phi \partial_\theta + 
   A_{22} \Theta \partial_\phi \eqno(22)$$
$$ \partial_\phi \Theta = B_{11} \Theta \partial_\phi + 
   B_{12} \Phi \partial_\theta $$
$$\partial_\phi \Phi = B_{21} \Phi \partial_\phi + 
  B_{22} \Theta \partial_\theta. $$
From the fact that 
$$\partial_i(\theta^j \Theta^k) = \delta^i_j \delta^k_l \Theta^k \eqno(23)$$
where $\partial_1 = \partial_\theta$, $\theta^1 = \theta$, 
$\Theta^1 = \Theta$, etc. we have 
$$A_{11} = 1 \qquad A_{22} = 0 \qquad B_{21} = 1 \qquad B_{12} = 0$$
$$F_{11} A_{21} + F_{12} A_{12} = 1 \qquad F_{21} A_{12} + F_{22} A_{21} = 0 $$
$$F_{21} B_{11} + F_{22} B_{22} = 1 \qquad F_{11} B_{22} + F_{12} B_{11} = 0. 
  \eqno(24)$$
Thus there are two types of solutions. 

{\it Type I} 
$$ \partial_\theta \Theta = \Theta \partial_\theta + 
   (1 - p^{-1}q^{-1}) \Phi \partial_\phi \qquad 
  \partial_\theta \Phi = q^{-1} \Phi \partial_\theta \eqno(25\mbox{a})$$
$$ \partial_\phi \Theta = p^{-1} \Theta \partial_\phi \qquad 
  \partial_\phi \Phi = \Phi \partial_\phi $$
and 

{\it Type II} 
$$ \partial_\theta \Theta = \Theta \partial_\theta \qquad 
   \partial_\theta \Phi = p \Phi \partial_\theta \eqno(25\mbox{b})$$
$$ \partial_\phi \Theta = q \Theta \partial_\phi \qquad 
  \partial_\phi \Phi = \Phi \partial_\phi + 
  (1 - pq) \Theta \partial_\theta. $$

We now should compute the R-matrix for these two types of solutions satisfying 
the Yang-Baxter equations 
$$R_{12} R_{13} R_{23} = R_{23} R_{13} R_{12} \qquad 
  \hat{R}_{12} \hat{R}_{23} \hat{R}_{12} = 
  \hat{R}_{23} \hat{R}_{12} \hat{R}_{23} \eqno(26)$$ 
where 
$$\hat{R}^{ij}_{kl} =  R^{ji}_{kl}. \eqno(27)$$
From the definition of $R$ matrix for the Yang-Baxter equation 
$$\theta^i \Theta^j = R^{ji}_{kl} \Theta^k \theta^l\eqno(28)$$
we have 
$$R 
 = \left(\matrix{ 1 &  0     &   0    & 0 \cr 
                  0 & F_{21} & F_{22} & 0 \cr
                  0 & F_{12} & F_{11} & 0 \cr 
                  0 &  0     & 0      & 1 \cr }\right) 
  = (R^{ij}{}_{kl}). \eqno(29)$$
The $R$ matrix for type I is 
$$R 
 = \left(\matrix{ 1 & 0 &   0    & 0 \cr 
                  0 & p & 1 - pq & 0 \cr
                  0 & 0 &   q    & 0 \cr 
                  0 & 0 &   0    & 1 \cr }\right) \eqno(30\mbox{a})$$
and for type II is 
$$R 
 = \left(\matrix{ 1 &  0               &   0    & 0 \cr 
                  0 & q^{-1}           & 0      & 0 \cr
                  0 & 1 - p^{-1}q^{-1} & p^{-1} & 0 \cr 
                  0 &  0               & 0      & 1 \cr }\right). 
   \eqno(30\mbox{b})$$
It can be checked that both matrices in (30) satisfy the Yang-Baxter 
equation (26). Note that the matrix in (30a) becomes the transposition of 
the inverse of the matrix in (30b) for $p = q$. 

We formulate the differential calculus, with the R-matrix, as follows: 

{\it Type I, II} 
$$\theta^i \theta^j = - {1\over {pq}} \hat{R}^{ij}{}_{kl} \theta^k \theta^l 
  \qquad 
  \theta^i \theta^j = - pq \hat{R}^{ij}{}_{kl} \theta^k \theta^l $$
$$\theta^i \Theta^j = \hat{R}^{ij}{}_{kl} \Theta^k \theta^l $$
$$\partial_i \theta^j = \delta^j{}_i - \hat{R}^{jk}{}_{il} \theta^l \partial_k 
  \eqno(31)$$
$$\partial_i \Theta^j = (\hat{R}^{-1})^{kj}{}_{li} \Theta^l \partial_k $$
$$\partial_i \partial_j = - {1\over {pq}} \hat{R}^{kl}{}_{ji} 
  \partial_l \partial_k \qquad 
\partial_i \partial_j = - pq \hat{R}^{kl}{}_{ji} \partial_l \partial_k. $$

To deform the group $GL(2)$, as two-parametric, we can use both R-matrices 
above, too. The equation 
$$\hat{R} T_1 T_2 = T_1 T_2 \hat{R} \eqno(32)$$ 
gives$^8$ the $(p,q)$-commutation relations between the matrix elements of any 
matrix $T$ in $GL_{p,q}(2)$, where $T_1 = T \otimes I$ and $T_2 = I \otimes T$ 
in the standard notation. One can verify that the relations (1), (2) and (11) 
are preserved under the transformations 
$$\theta \longrightarrow a \theta + b \phi \qquad 
  \Theta \longrightarrow a \Theta + b \Phi $$
$$\phi \longrightarrow c \theta + d \phi \qquad 
  \Phi \longrightarrow c \Theta + d \Phi \eqno(33)$$
provided (32) satisfied. Here, it assummed that the coordinates and 
differentials commute with the matrix elements $a$, $b$, $c$ and $d$. In the 
same sense the relations (14) [or (16)] are preserved under the action of 
the quantum matrix $T \in GL_{p,q}(2)$, and so. Consequently, 
two-parameter calculus is covariant under the action of the quantum 
group $GL_{p,q}(2)$. 

We finally obtain a two-parameter deformation of fermionic oscillator algebra. 
For this, we assume that 
$$p = \bar{q}. \eqno(34)$$
We then introduce two fermionic oscillators, $B_1$ and $B_2$, and make 
identification 
$$\theta \longrightarrow B_1^+ \qquad \partial_\theta \longrightarrow B_1 $$
$$\phi \longrightarrow B_2^+ \qquad \partial_\phi \longrightarrow B_2 
  \eqno(35)$$
we construct the $(p,q)$-deformed fermionic oscillator algebra as follows: 
$$B_1 B_2 + q B_2 B_1 = 0 \qquad B_1^2 = 0 = {B_1^+}^2 $$
$$B_1^+ B_2^+ + p^{-1} B_2^+ B_1^+ = 0 \qquad B_2^2 = 0 = {B_2^+}^2 \eqno(36)$$

{\it Type I}
$$B_1 B_2^+ + p B_2^+ B_1 = 0 \qquad B_2 B_1^+ + q B_1^+ B_2 = 0 $$
$$B_1 B_1^+ + B_1^+ B_1 = 1 \qquad 
  B_2 B_2^+ + B_2^+ B_2 = 1 + (pq - 1) B_1^+ B_1 \eqno(37\mbox{a})$$

{\it Type II}
$$B_1 B_2^+ + q^{-1} B_2^+ B_1 = 0 \qquad B_2 B_1^+ + p^{-1} B_1^+ B_2 = 0 $$
$$B_1 B_1^+ + B_1^+ B_1 = 1 + (p^{-1}q^{-1} - 1) B_2^+ B_2 \qquad 
  B_2 B_2^+ + B_2^+ B_2 = 1. \eqno(37\mbox{b})$$

{\bf Acknowledgment}

This work was supported in part by {\bf T. B. T. A. K.} the 
Turkish Scientific and Technical Research Council.

\end{document}